\documentclass[twoside,a4paper,12pt,reqno] {amsart}
\usepackage{amsmath,mathrsfs}
\usepackage{epsfig}\usepackage{psfrag}\usepackage{subfigure}
\usepackage{tabularx}\usepackage{rotating}
\newtheorem{theorem}{Theorem}

\newtheorem{lemma}{Lemma}

\def\cal{\mathcal}

\begin{document}
\title[The weak Banach-Saks Property of the Space $(L_\mu^p)^m$]
{The weak Banach-Saks Property of the Space $(L_\mu^p)^m$}

\author{Zhenglu Jiang}
\address{Department of Mathematics, Zhongshan University, 
Guangzhou 510275, China}
\email{mcsjzl@mail.sysu.edu.cn}

\author{Xiaoyong Fu}
\address{Department of Mathematics, Zhongshan University, 
Guangzhou 510275, China
}
\email{mcsfxy@mail.sysu.edu.cn}


\baselineskip 16pt

\date{\today.}


\keywords{convex functions; measure spaces; Lebesgue integrals; Banach spaces}

\begin{abstract}
In this paper we show the weak Banach-Saks property of  
the Banach vector space $(L_\mu^p)^m$ generated by $m$ $L_\mu^p$-spaces for $1\leq p<+\infty,$ 
where $m$ is any given natural number.  
When $m=1,$ this is the famous Banach-Saks-Szlenk theorem.  
By use of this property,  we also present inequalities for integrals of functions   
that are the composition of nonnegative continuous convex functions 
on a convex set of a vector space ${\bf R}^m$ 
and vector-valued functions in a weakly compact subset of
the space $(L_\mu^p)^m$ for $1\leq p<+\infty$  
and inequalities when 
these vector-valued functions are in a  weakly* compact subset of   
the product space $(L_\mu^\infty)^m$ generated by $m$ $L_\mu^\infty$-spaces. 
\end{abstract}

\maketitle

\section{Introduction}
\label{intro} 
We begin with some notations and definitions used throughout this paper. $m$ and $n$ are natural numbers, 
${\bf R}$ denotes the real number system, ${\bf R}^n$ is the usual vector 
space of real $n$-tuples $x=(x_1,x_2,\cdots,x_n),$ 
$\mu$  is a nonnegative Lebesgue measure of ${\bf R}^n;$ 
$L_\mu^p({\bf R}^n)$ represents a Banach space of any measurable function $\hat{u}=\hat{u}(x)$ 
with its finite norm 
\begin{equation}
{\parallel}\hat{u}{\parallel_p}=\left(\int_{{\bf R}^n}|\hat{u}(x)|^p d\mu\right)^{\frac{1}{p}}
\nonumber\label{norm}
\end{equation} 
for any given $p\in [1,+\infty),$  and $(L_\mu^p({\bf R}^n))^m$ denotes a Banach vector space 
where each measurable vector-valued function $u=u(x)$ with $m$ components 
$\hat{u}^{(j)}=\hat{u}^{(j)}(x)$ ($j=1,2,\cdots,m$) in $L_\mu^p({\bf R}^n)$ is given by  $u=(\hat{u}^{(1)},\hat{u}^{(2)},\cdots,\hat{u}^{(m)})$ 
and its norm is defined by ${\parallel}u{\parallel_p}=(\sum\limits_{j=1}^m{\parallel}\hat{u}^{(j)}{\parallel_p^p})^{1/p};$   
similarly, $L_\mu^\infty({\bf R}^n)$ represents a Banach space of any measurable function $\hat{u}=\hat{u}(x)$ 
with  its finite norm 
\begin{equation}
{\parallel}\hat{u}{\parallel_\infty}= 
\mathop{\hbox{ess}\sup}\limits_{x\in {\bf R}^n} |\hat{u}(x)| 
\hbox{  } (\hbox{or say  }{\parallel}\hat{u}{\parallel_\infty}  
=\mathop{\inf\limits_ {E\subseteq {\bf R}^n}}
\limits_{\mu(E^C)=0}\sup\limits_{x\in E} |\hat{u}(x)| ) 
\nonumber\label{norminf}
\end{equation} 
where $E^C$ represents the complement set of $E$ in ${\bf R}^n,$ 
 and  $(L_\mu^\infty({\bf R}^n))^m$ denotes a Banach vector space 
where each measurable vector-valued function  $u=u(x)$ with $m$ components 
$\hat{u}_j=\hat{u}_j(x)$ ($i=1,2,\cdots,m$) in $L_\mu^\infty({\bf R}^n)$ is denoted by  $u=(\hat{u}^{(1)},\hat{u}^{(2)},\cdots,\hat{u}^{(m)})$     
and its norm is defined by ${\parallel}u{\parallel_\infty}=\sum\limits_{j=1}^m{\parallel}\hat{u}^{(j)}{\parallel_\infty}.$ 

If a function $u$ and a sequence $\{u_i\}_{i=1}^{+\infty}$ in $(L_\mu^p({\bf R}^n))^m$ 
are assumed to satisfy the fact that $\lim\limits_{i\to\infty}{\parallel}u_i-u{\parallel_p}=0,$ then 
this sequence $\{u_i\}_{i=1}^{+\infty}$ is said to 
be strongly convergent in  
 $(L_\mu^p({\bf R}^n))^m$ 
to $u.$  
Similarly, if a function $u$ and  
a sequence $\{u_i\}_{i=1}^{+\infty}$ in $(L_\mu^\infty({\bf R}^n))^m$ are assumed to
have the property that $\lim\limits_{i\to\infty}{\parallel}u_i-u{\parallel_\infty}=0,$
then 
this sequence $\{u_i\}_{i=1}^{+\infty}$ is said to 
be strongly convergent in  $(L_\mu^\infty({\bf R}^n))^m$ 
to $u.$      
 
Besides the convergence given above, we consider a definition of weak  
convergence of a sequence in  $(L_\mu^p({\bf R}^n))^m.$  
  Assume that   
$q=p/(p-1) $  as  $p\in (1, +\infty) $ and that  
$q=\infty$ as  $p=1.$ 
If a function $\hat{u}$ and
 a sequence $\{\hat{u}_i\}_{i=1}^{+\infty}$ in $L_\mu^p({\bf R}^n)$ 
 have the following relation:     
\begin{equation}
 \mathop{\lim}\limits_{i\rightarrow+\infty}
\int_{{\bf R}^n}\hat{u}_i\hat{v}d\mu=\int_{{\bf R}^n}\hat{u}\hat{v}d\mu
\label{weqp}
\end{equation}
for all $\hat{v}\in L_\mu^q({\bf R}^n),$ then the sequence $\{\hat{u}_i\}_{i=1}^{+\infty}$ 
is said to be weakly convergent in 
$L_\mu^p({\bf R}^n)$ to $\hat{u}.$   
If $\{\hat{u}_{i}^{(j)}\}_{i=1}^{+\infty}$ 
is weakly convergent in $L_\mu^p({\bf R}^n)$ to $\hat{u}^{(j)}$ 
 for all $j=1,2,\cdots,m$ as $i$ goes to $\infty,$  
then a sequence $\{u_i=(\hat{u}_{i}^{(1)},\hat{u}_{i}^{(2)},\cdots,\hat{u}_{i}^{(m)})\}_{i=1}^{+\infty}$ is said to 
be weakly convergent in  
 $(L_\mu^p({\bf R}^n))^m$ 
to $u=(\hat{u}^{(1)},\hat{u}^{(2)},\cdots,\hat{u}^{(m)}).$  

Similarly, we introduce a definition of weak*   
convergence of a sequence in $L_\mu^\infty({\bf R}^n).$  
If a function $\hat{u}$ and
 a sequence $\{\hat{u}_i\}_{i=1}^{+\infty}$  in $L_\mu^\infty({\bf R}^n)$  
 satisfy  the equality (\ref{weqp}) 
for all $\hat{v}\in L_\mu^1({\bf R}^n),$ then 
the sequence $\{\hat{u}_i\}_{i=1}^{+\infty}$ is said to be weakly* convergent in 
$L_\mu^\infty({\bf R}^n)$ to $\hat{u}.$  
If $\{\hat{u}_{i}^{(j)}\}_{i=1}^{+\infty}$ 
is weakly* convergent in $L_\mu^\infty({\bf R}^n)$ to $\hat{u}^{(j)}$ 
 for all $j=1,2,\cdots,m$ as $i$ goes to $\infty,$  
then a sequence $\{u_i=(\hat{u}_{i}^{(1)},\hat{u}_{i}^{(2)},\cdots,\hat{u}_{i}^{(m)})\}_{i=1}^{+\infty}$ is said to 
be weakly* convergent in  
 $(L_\mu^\infty({\bf R}^n))^m$ 
to $u=(\hat{u}^{(1)},\hat{u}^{(2)},\cdots,\hat{u}^{(m)}).$  

Assume that $Y$ is a Banach space. 
$Y$ is said to be of the weak Banach-Saks property 
if any sequence $\{y_i\}_{i=1}^{+\infty}$ weakly convergent  in $Y$ to $y$ contains 
a subsequence $\{y_{i_r}\}_{r=1}^{+\infty}$ 
such that $(\sum_{r=1}^ky_{i_r})/k$ converges strongly in $Y$ to $y.$

Banach and Saks \cite{bs} first proved that  $L_\mu^p({\bf R}^n)$ 
has the weak Banach-Saks property for the $1<p<+\infty$ case in 1930 
and then the similar result for $L_\mu^1({\bf R}^n)$ was showed by 
Szlenk \cite{ws} in 1965. This result about $L_\mu^p({\bf R}^n)$ 
is the famous Banach-Saks-Szlenk theorem. 
Now there is not yet this result about the 
Banach vector space $(L_\mu^p({\bf R}^n))^m$ when $m\not=1.$ 
The aim of this paper is to extend the Banach-Saks-Szlenk theorem to the case of 
the vector space $(L_\mu^p({\bf R}^n))^m$ generated by $m$ $L_\mu^p({\bf R}^n)$-spaces for $1\leq p<+\infty$ 
and show that $(L_\mu^p({\bf R}^n))^m$ has the weak Banach-Saks property for any fixed natural number $m.$ 
An application of this property is to show inequalities for integrals of functions   
that are the composition of nonnegative continuous convex functions 
on a convex set of a vector space ${\bf R}^m$ 
and vector-valued functions in a weakly compact subset of  
 a Banach vector space generated by $m$ $L_\mu^p$-spaces for $1\leq p<+\infty$ 
and  inequalities when these vector-valued functions are in a  weakly* compact subset of 
 a Banach vector space generated by $m$ $L_\mu^\infty$-spaces. 

\section{The Weak Banach-Saks Property}
\label{bsp} 
 A detail description of the weak Banach-Saks property of $(L_\mu^p({\bf R}^n))^m$ for any fixed natural number $m$ 
is as follows: 
\begin{theorem}
Given a real number $p$ in $ [1,+\infty).$  
Assume that a sequence $\{u_i=u_i(x)\}_{i=1}^{+\infty}$ converges weakly in $(L_\mu^p({\bf R}^n))^m$ to 
$u=u(x).$  Then this sequence contains a subsequence $\{u_{i_r}\}_{r=1}^{+\infty}$ 
 with its arithmetic means $\frac{1}{k}\sum_{r=1}^ku_{i_r}$ strongly 
convergent  in $(L_\mu^p({\bf R}^n))^m$ to $u$  as $k$ goes to infinity. 
\label{bss}
\end{theorem}
We can below show Theorem~\ref{bss} using the two following techniques with only minor adjustments: 
one is given by Banach and Saks \cite{bs} for any fixed $p\in(1,+\infty),$ 
and another by Szlenk \cite{ws} in the case when $p=1.$  

To prove Theorem \ref{bss} for any fixed $p\in(1,+\infty),$ 
we have first to introduce the following lemma: 
\begin{lemma}[\cite{bs}]
Let $a$ and $b$ be any real numbers  and $1<p<+\infty.$ Then
\begin{equation} 
|a+b|^p\leq |a|^p+p|a|^{p-1}[sgn(a)]b+A|b|^p+B(p,a,b).
\label{ineq01f}
\end{equation}
Here, $A$ is a positive constant independent of $a$ and $b;$ $sgn(\tau)$ is defined as follows: 
$sgn(0)=0,$ $sgn(\tau)=1$ as $\tau>0$ and $sgn(\tau)=-1$ as $\tau<0;$
$B(p,a,b)=0$ as $p\in (1,2]$ and $B(p,a,b)=\sum_{i=2}^{E(p)}{p \choose i}|a|^{p-i}|b|^i$ as $p\in (2,+\infty),$ 
where $E(p)$ is the largest natural number less than $p.$ 
\label{lem01f}
\end{lemma}
The proof of this lemma can be found in \cite{bs}. Using the inequality (\ref{ineq01f}), we can get a similar result 
to that given by Banach and Saks \cite{bs}. This result is as follows:
\begin{lemma}
Assume that $p>1$ and that a sequence $\{\hat{u}_i=\hat{u}_i(x)\}_{i=1}^{+\infty}$ in $L_\mu^p({\bf R}^n)$ satisfies 
\begin{equation} 
\int_{{\bf R}^n}|\hat{u}_i(x)|^pd\mu\leq 1
\label{ineq02f}
\end{equation}
for all $i\geq 1.$ Put $\hat{s}_k(x)=\sum_{i=1}^k\hat{u}_i(x).$ Then
\begin{equation} 
\int_{{\bf R}^n}|\hat{s}_k(x)|^pd\mu\leq C(k)+p\int_{{\bf R}^n}|\hat{s}_{k-1}(x)|^{p-1}[sgn(\hat{s}_{k-1}(x))]\hat{u}_k(x)d\mu
+\int_{{\bf R}^n}|\hat{s}_{k-1}(x)|^p d\mu, 
\label{ineq03f}
\end{equation}
where, $C(k)=A+Bk^{p-2},$ $A$ and $B$ are positive constants independent of $k$ and 
$\{\hat{u}_i\}_{i=1}^{+\infty},$  $sgn(\tau)$ is defined as in Lemma \ref{lem01f}. 
\label{lem02f}
\end{lemma}
\begin{proof}
Insert $a=\hat{s}_{k-1}(x)$ and $b=\hat{u}_k(x)$ into the inequality (\ref{ineq01f}) 
and integrate all its terms over the space ${\bf R}^n.$ 
Then, by (\ref{ineq02f}), we can know that (\ref{ineq03f}) holds in the $1<p\leq 2$ case 
and that 
\begin{eqnarray} 
\int_{{\bf R}^n}|\hat{s}_k(x)|^pd\mu\leq \int_{{\bf R}^n}|\hat{s}_{k-1}(x)|^p d\mu
+p\int_{{\bf R}^n}|\hat{s}_{k-1}(x)|^{p-1}[sgn(\hat{s}_{k-1}(x))]\hat{u}_k(x)d\mu \nonumber \\
+\sum\limits_{i=2}^{E(p)}\left( _i^p\right)\int_{{\bf R}^n}|\hat{s}_{k-1}(x)|^{p-i}|\hat{u}_k(x)|^i d\mu+A 
\label{ineq04f}
\end{eqnarray}
for all $p>2.$ 
Notice that $\int_{{\bf R}^n}|\hat{s}_{k-1}(x)|^{p-i}|\hat{u}_k(x)|^i d\mu\leq k^{p-i}$ for all $i\leq p;$ 
this can be obtained by first using the H\"{o}lder inequality and then the Minkowski one 
with the help of the condition (\ref{ineq02f}). Take $B=\sum_{i=2}^{E(p)}{p \choose i}.$  
It can be then found that (\ref{ineq04f}) gives (\ref{ineq03f}) for $p> 2.$ 
This hence completes our proof. 
\end{proof} 

We below give the proof of Theorem \ref{bss} for any fixed $p\in (1,+\infty).$ 
To do this, it suffices to consider the case when $m=2.$  Let us first denote 
all the vector-valued functions $u_i$ of this sequence in $(L_\mu^p({\bf R}^n))^2$
by $u_i=(\hat{u}_{i}^{(1)},\hat{u}_{i}^{(2)}),$ 
where $\hat{u}_{i}^{(1)}=\hat{u}_{i}^{(1)}(x)$ and $\hat{u}_{i}^{(2)}=\hat{u}_{i}^{(2)}(x)$ 
represent two functions in $L_\mu^p({\bf R}^n)$  for any natural number $i.$ 
Since any weak convergent sequence in $L_\mu^p({\bf R}^n)$ is bounded,   
we may first assume without loss of generality that 
all the functions $u_i$ of the considered sequence satisfy 
\begin{equation} 
{\parallel}u_i{\parallel_p^p}\leq 1
\label{ineq05f}
\end{equation}
for all $i\geq 1.$ We may also assume without loss of generality that 
 this sequence $\{u_i\}_{i=1}^{+\infty}$ converges weakly in $(L_\mu^p({\bf R}^n))^2$ to zero. 
Then, by recursion, we can determine a subsequence 
$\{u_{i_r}=(\hat{u}_{i_r}^{(1)},\hat{u}_{i_r}^{(2)})\}_{i=1}^{+\infty}$ ($i_1=1$).   
This recursive process can be roughly divided into two steps and they are described as follows.
The first step to do this  is to take $\hat{s}_k^{(j)}(x)=\sum_{r=1}^k\hat{u}_{i_r}^{(j)}(x)$ for $j=1,2$ under 
 the assumption that $k$ previous terms  
$\{u_{i_r}=(\hat{u}_{i_r}^{(1)},\hat{u}_{i_r}^{(2)})\}_{i=1}^{k}$ of this subsequence is determined.  
 It can be then known that $\hat{s}_k^{(j)}(x)\in L_\mu^p({\bf R}^n)$ and that   
$|\hat{s}_{k}^{(j)}(x)|^{p-1}[sgn(\hat{s}_{k}^{(j)}(x))]\in L_\mu^{p/(p-1)}({\bf R}^n)$ for $j=1,2;$
since these functions $\hat{u}_{i}^{(j)}$ converge weakly in $ L_\mu^p({\bf R}^n)$ to zero for $j=1,2,$  
there exists a natural number $i_k$ such that 
\begin{equation} 
\int_{{\bf R}^n}|\hat{s}_{k}^{(j)}(x)|^{p-1}[sgn(\hat{s}_{k}^{(j)}(x))]\hat{u}_i^{(j)}(x)d\mu\leq 1
\label{ineq06f}
\end{equation}
for all $i>i_k$ and $j=1,2;$ thus the second one is to define the subscript $i_{k+1}$ of the next term 
to be one of
all the natural numbers $i$ satisfying the condition given by (\ref{ineq06f}). 

Then, by  (\ref{ineq06f}), we can know that for all $k>1$ and $j=1,2,$
\begin{equation} 
\int_{{\bf R}^n}|\hat{s}_{k-1}^{(j)}(x)|^{p-1}[sgn(\hat{s}_{k-1}^{(j)}(x))]\hat{u}_{i_k}^{(j)}(x)d\mu\leq 1.
\label{ineq07f}
\end{equation}
Combining (\ref{ineq05f}) and (\ref{ineq07f}) and using Lemma \ref{lem02f}, 
we can show that for $j=1,2,$  
\begin{equation} 
\int_{{\bf R}^n}|\hat{s}_{k}^{(j)}(x)|^{p}d\mu\leq (A+p)k+Bk^{p-2}+1, 
\nonumber
\end{equation}
thus giving  
\begin{equation} 
\lim\limits_{k\to\infty}\int_{{\bf R}^n}\left|\frac{\hat{s}_{k}^{(j)}(x)}{k}\right|^{p}d\mu=0 \hbox{ for } j=1,2.
\nonumber
\end{equation}
This hence completes our proof of Theorem \ref{bss} for any given $p\in (1,+\infty).$  

Now it remains to prove Theorem \ref{bss} when $p=1.$ 
To do this, we first recall a lemma as follows:
\begin{lemma}[\cite{ws}] 
Assume that $\hat{u}_i$ belongs to the Banach space $ L[0,1]$ for any natural number $i$ and 
converges weakly to zero as $i$ goes to infinity.  
Then, for any given $\varepsilon>0,$ there exists a sequence of indexes $i_r$ such that 
$
\mathop{\overline{\lim}}\limits_{k\to\infty}\mathop{\sup}\limits_{i_1<\cdots<i_k}
\frac{1}{k}{\parallel}\sum\limits_{r=1}^k\hat{u}_{i_r}{\parallel_{L}}\leq \varepsilon
$
where ${\parallel}\cdot{\parallel_{L}}$ represents the norm of the Banach space $L[0,1].$ 
\label{lem01s}
\end{lemma}
Lemma \ref{lem01s} and its proof were shown by Szlenk \cite{ws} in 1965. 
By using a similar proof to that given by Szlenk, it can be found that 
Lemma \ref{lem01s} still holds if 
 $ L[0,1]$ is replaced by the Banach space $ L([0,1]^n).$ 
Since $L({\bf R}^n)$ is isometric to $ L([0,1]^n)$ (see \cite{pw}, Page 83), 
we can easily deduce that 
\begin{lemma} 
Assume that $\hat{u}_i$ belongs to the Banach space $L({\bf R}^n)$ for any natural number $i$ and 
converges weakly to zero as $i$ goes to infinity.   
Then, for any given $\varepsilon>0,$ there exists a sequence of indexes $i_r$ such that 
$ 
\mathop{\overline{\lim}}\limits_{k\to\infty}\mathop{\sup}\limits_{i_1<\cdots<i_k}
\frac{1}{k}{\parallel}\sum\limits_{r=1}^k\hat{u}_{i_r}{\parallel_{1}}\leq \varepsilon.
$
\label{lem02s}
\end{lemma}
Using Lemma \ref{lem02s}, we can get the following result: 
\begin{lemma} 
Assume that $u_i$ belongs to the Banach space $(L({\bf R}^n))^m$ for any natural number $i$ and 
converges weakly to zero as $i$ goes to infinity.   
Then, for any given $\varepsilon>0,$ there exists a sequence of indexes $i_r$ such that 
\begin{equation} 
\mathop{\overline{\lim}}\limits_{k\to\infty}\mathop{\sup}\limits_{i_1<\cdots<i_k}
\frac{1}{k}{\parallel}\sum\limits_{r=1}^ku_{i_r}{\parallel_{1}}\leq \varepsilon.
\end{equation}
\label{lem03s}
\end{lemma}

We can below prove Theorem \ref{bss} when $p=1.$ 
To do this, it suffices to consider the case of $u=0.$ By Lemma \ref{lem03s}, for any given $l\geq 1,$ 
there exists a sequence of indexes $i_{l,r}$ such that 
\begin{equation} 
\mathop{\overline{\lim}}\limits_{k\to\infty}\mathop{\sup}\limits_{s_1<\cdots<s_k}
\frac{1}{k}{\parallel}\sum\limits_{r=1}^ku_{i_{l,s_r}}{\parallel_{1}}\leq \frac{1}{l}.
\label{ineq01s}
\end{equation}
Assume that the sequence of indexes $i_{l+1,r}$ 
is a subsequence of the sequence of indexes $i_{l,r}.$ 
Denote by $\{u_{i_{r}}\}_{r=1}^{+\infty}$ a sequence of indexes $i_{r}=i_{r,r}$ corresponding to the condition (\ref{ineq01s}). 
Then we can know that this sequence $\{u_{i_{r}}\}_{r=1}^{+\infty}$ satisfies 
\begin{equation}
\frac{1}{k}{\parallel}\sum\limits_{r=1}^ku_{i_r}{\parallel_{1}}\leq 
\frac{1}{k}{\parallel}\sum\limits_{r=1}^l u_{i_r}{\parallel_{1}}
+\frac{1}{k-l}{\parallel}\sum\limits_{r=1}^{k-l}u_{i_{l+r}}{\parallel_{1}}
\nonumber
\end{equation} 
for all $k>l.$ It follows that 
\begin{equation}
\mathop{\overline{\lim}}\limits_{k\to\infty}\frac{1}{k}{\parallel}\sum\limits_{r=1}^ku_{i_r}{\parallel_{1}}\leq
\mathop{\overline{\lim}}\limits_{k\to\infty}\frac{1}{k-l}{\parallel}\sum\limits_{r=1}^{k-l}u_{i_{l+r}}{\parallel_{1}}
=\mathop{\overline{\lim}}\limits_{k\to\infty}\frac{1}{k}{\parallel}\sum\limits_{r=1}^{k}u_{i_{l+r}}{\parallel_{1}}.
\label{ineq02s}
\end{equation} 
Since $\{u_{i_{r}}\}_{r=1}^{+\infty}$ is a subsequence of this sequence $\{u_{i_{l,r}}\}_{r=1}^{+\infty},$ 
we have 
\begin{equation}
\frac{1}{k}{\parallel}\sum\limits_{r=1}^{k}u_{i_{l+r}}{\parallel_{1}}\leq
\mathop{\sup}\limits_{s_1<\cdots<s_k}
\frac{1}{k}{\parallel}\sum\limits_{r=1}^ku_{i_{l,s_r}}{\parallel_{1}}.
\label{ineq03s}
\end{equation} 
Combining (\ref{ineq01s}), (\ref{ineq02s}) and (\ref{ineq03s}), we can know that
\begin{equation}
\mathop{\overline{\lim}}\limits_{k\to\infty}\frac{1}{k}{\parallel}\sum\limits_{r=1}^{k}u_{i_{r}}{\parallel_{1}}\leq\frac{1}{l}
\hbox{ for } l=1,2,\cdots.
\nonumber
\end{equation} 
This implies that $\frac{1}{k}\sum_{r=1}^{k}u_{i_{r}}$ 
converges strongly in the Banach space $(L({\bf R}^n))^m$ to zero 
as $k$ tends to infinity. 
Our proof is hence finished.

We can also extend Theorem \ref{bss} to a more general case, that is,  
\begin{theorem}
Given a measure space $(X,{\cal A},\mu)$ and a real number $p$ in $ [1,+\infty).$ 
Assume that   
 $\{\hat{u}_i^{(j)}\}_{i=1}^{+\infty}$ converges weakly in $(L_\mu^p(X))^m$ to 
$\hat{u}^{(j)}$ for $j=1,2,\cdots,m.$ 
Take $u=(\hat{u}^{(1)},\hat{u}^{(2)},\cdots,\hat{u}^{(m)})$ and 
$u_i=(\hat{u}_i^{(1)},\hat{u}_i^{(2)},\cdots,\hat{u}_i^{(m)})$ for any natural number $i.$
Then this sequence $\{u_i\}_{i=1}^{+\infty}$ contains a subsequence $\{u_{i_r}\}_{r=1}^{+\infty}$ 
 with its arithmetic means $\frac{1}{k}\sum_{r=1}^ku_{i_r}$ strongly 
convergent  in $(L_\mu^p(X))^m$ to 
$u$ as $k$ goes to infinity. 
\label{bssg}
\end{theorem}

When $m=1,$ this result appears in the book 
of Benedetto \cite{bjj} and 
there is an explanation of its proof, that is, it is the same as 
given by Banach and Saks for any $p\in (1,+\infty)$  
and by Szlenk in the case when $p=1.$ 
We below give a simple proof of Theorem \ref{bssg}. 
Notice that the separable space $L_\mu^p(X)$  
is isometric to $L_\mu^p[0,1]$ for any $p\in [1,+\infty)$ (see \cite{pw}).  
It is then clear that $(L_\mu^p(X))^m$ 
is also isometric to $(L_\mu^p[0,1])^m$ for any fixed natural number $m.$ 
By Theorem \ref{bss}, Theorem \ref{bssg} thus follows.

\section{Application to Inequalities for Integrals}
\label{aii}  
By use of  the weak Banach-Saks property of $(L_\mu^p)^m$, we can show inequalities for integrals of functions   
which are the composition of nonnegative  continuous convex functions 
on a convex set of a vector space ${\bf R}^m$ 
and vector-valued functions in a weakly compact subset  of 
 a Banach vector space generated by $m$ $L_\mu^p$-spaces for any given $p\in [1,+\infty).$ 
That is the following
\begin{theorem}  
Suppose that a sequence $\{u_i\}_{i=1}^{+\infty}$ 
 weakly converges in 
 $(L_\mu^p({\bf R}^n))^m$ to $u$   
 as $i$ goes to infinity,  where $p\in [1,+\infty)$ and $m$ and $n$ are two positive integers.           
Assume that all the values of $u$ and $u_i$ ($i=1,2,3,\cdots$)  
belong to a convex set $K$ in ${\bf R}^m$  and that 
$f(w)$ is a nonnegative continuous convex function from $K$ to  ${\bf R}.$ 
 Then 
\begin{equation}
\mathop{\underline{\lim}}\limits_{i\rightarrow+\infty}
 \int_{\Omega}f(u_i)d\mu\geq\int_{\Omega}f(u)d\mu \label{ineq}
\end{equation}
for any measurable set $\Omega\subseteq {\bf R}^n.$
\label{th01}
\end{theorem}

The estimates of integrals of this kind of composite function  
is interesting and important in many application areas such as  
the existence of solutions of differential equations 
(e.g., see \cite{dl} and \cite{yl}). 
A similar result  was shown by Jiang {et.}~al \cite{jipam} 
if  $K$ is assumed to be an open convex set of  ${\bf R}^m$ instead; 
 Egorov's theorem is used into their proof 
except for the weak Banach-Saks property of $(L_\mu^p)^m.$   
Meanwhile, a simple proof of another similar one was given in \cite{jmaa} 
when $K$ is set to be ${\bf R}^m;$ this proof requires  
 the weak Banach-Saks property of $(L_\mu^p)^m$ but it does not give 
any proof of this property; 
it only shows  the case when $m=1$ in Theorem~\ref{bssg}. 
The former device is valid for only an open convex set $K$ while the latter one   
 is suitable to show inequalities for integrals of these composite functions 
in a more general case, or more precisely speaking, 
this case is for any convex set $K$ in ${\bf R}^m.$ 
Therefore it is still very necessary to show Theorem~\ref{th01} and its proof. 

It is worth mentioning that some properties of convex functions and weakly compact sets 
 can be found in the literature (e.g., see \cite{bjj}, \cite{srl}, \cite{rrt}, \cite{pw} and  \cite{yk}). 

\begin{proof}[Proof of Theorem \ref{th01}]
Put $\alpha_i=\int_{\Omega}f(u_i)d\mu$ ($i=1,2,\cdots $) and 
$\alpha=\mathop{\underline{\lim}}\limits_{i\rightarrow +\infty}\int_{\Omega}f(u_i)d\mu$ 
for all $\Omega\subseteq {\bf R}^n.$  
Then there exists a subsequence of $\{\alpha_i\}_{i=1}^{+\infty}$ such that 
this subsequence, denoted without loss of generality 
by $\{\alpha_i\}_{i=1}^{+\infty},$ converges to $\alpha$ as $i\rightarrow +\infty.$  

Since $u_i$ converges weakly in $(L_\mu^p({\bf R}^n))^m$ to $u$ for $1\leq p<+\infty,$ 
by Theorem~\ref{bss},  
it is easy to see that there exists a subsequence $\{u_{i_j}:j=1,2,\cdots\}$ 
such that $\frac{1}{k}\sum_{j=1}^ku_{i_j} \rightarrow u$ 
in $(L_\mu^p({\bf R}^n))^m$ for $1\leq p<+\infty$ as $k\rightarrow +\infty.$  
Thus there exists a subsequence of $\{\frac{1}{k}\sum_{j=1}^ku_{i_j}: k=1,2,\cdots\}$ 
such that this subsequence (also denoted without loss of generality by 
$\{\frac{1}{k}\sum_{j=1}^ku_{i_j}: k=1,2,\cdots\}$) satisfies that, as $k\rightarrow +\infty,$
\begin{equation}
\frac{1}{k}\sum_{j=1}^ku_{i_j} \rightarrow u  \hbox{  a.e. in  }{\bf R}^n.
\label{conv}
\end{equation} 
On the other hand, since all the values of $\{u_i\}_{i=1}^{+\infty}$ and $u$ 
belong to the convex set $K$  in ${\bf R}^m$ and 
$f(w)$ is a nonnegative continuous convex function from $K$ to  ${\bf R},$  we have
\begin{equation}
f(\frac{1}{k}\sum_{j=1}^ku_{i_j})\leq \frac{1}{k}\sum_{j=1}^kf(u_{i_j}). 
\label{conf1}
\end{equation}
By (\ref{conf1}) and Fatou's lemma, it follows that 
\begin{equation}
 \int_{\Omega}\mathop{\underline{\lim}}\limits_{k\rightarrow+\infty}
f(\frac{1}{k}\sum_{j=1}^ku_{i_j})d\mu
\leq\mathop{\underline{\lim}}\limits_{k\rightarrow+\infty}
\frac{1}{k}\sum_{j=1}^k\int_{\Omega}f(u_{i_j})d\mu. 
\label{conf2}
\end{equation}
Combining (\ref{conv}) and (\ref{conf2}), we can know that 
\begin{equation}
 \int_{\Omega}f(u)d\mu
\leq\mathop{\underline{\lim}}\limits_{k\rightarrow+\infty}
\frac{1}{k}\sum_{j=1}^k\int_{\Omega}f(u_{i_j})d\mu\equiv 
\mathop{\underline{\lim}}\limits_{k\rightarrow+\infty}
\frac{1}{k}\sum_{j=1}^k\alpha_{i_j}.
\label{conf3}
\end{equation}
Finally, by using the property of the convergence of 
$\alpha_i$ to $\alpha,$  
(\ref{conf3}) gives (\ref{ineq}). 
This completes our proof.
\end{proof}

Furthermore, we can give the following 
similar result for weakly* convergent sequences 
in $(L_\mu^\infty({\bf R}^n))^m:$   
\begin{theorem} 
Assume that a sequence $\{u_i\}_{i=1}^{+\infty}$ 
 weakly* converges in 
$(L_\mu^\infty({\bf R}^n))^m$ to $u$ as $i$ goes to infinity,  
where $m$ and $n$ are two positive integers. 
Assume that all the values of $u$ and $u_i$ ($i=1,2,3,\cdots$)   
belong to a convex set $K$ of  ${\bf R}^m$  and that 
$f(w)$ is a nonnegative continuous convex function from $K$ to  ${\bf R}.$ 
Then the inequality (\ref{ineq}) holds 
for any measurable set $\Omega\subseteq {\bf R}^n.$
\label{th02}
\end{theorem}
\begin{proof}  
Put $\Omega_R=\Omega\cap \{x: |x|<R,x\in {\bf R}^n\}.$ 
Then $\Omega_R$  is a bounded measurable set 
in ${\bf R}^n$ for all the fixed positive real number $R.$ 
Since $u_i\rightarrow{u}$ weakly* in 
$(L_\mu^\infty({\bf R}^n))^m,$ $u_i\rightarrow{u}$ weakly* in 
$(L_\mu^\infty(\Omega_R))^m.$ Hence, 
by $L^\infty(\Omega_R)\subset L^1(\Omega_R),$ 
it can be easily known that 
$u_i\rightarrow{u}$ weakly in 
$(L_\mu^1(\Omega_R))^m.$ 
Then, using the process of the proof of Theorem~\ref{th01},
we can get 
\begin{equation}
\mathop{\underline{\lim}}\limits_{i\rightarrow+\infty}
 \int_{\Omega_R}f(u_i)d\mu\geq\int_{\Omega_R}f(u)d\mu. \label{ineq01}
\end{equation}
It follows from the nonnegativity of the convex function $f$ that 
\begin{equation}
\mathop{\underline{\lim}}\limits_{i\rightarrow+\infty}
 \int_{\Omega}f(u_i)d\mu\geq\int_{\Omega_R}f(u)d\mu. \label{ineq02}
\end{equation}
Finally, by Lebesgue monotonous convergence theorem, 
as $R\rightarrow +\infty,$ (\ref{ineq02}) 
implies (\ref{ineq}).  Our proof is completed. 
\end{proof}

Also, removing the nonnegativity of $f(w)$ 
and assuming that the convex set $K$ is closed, by Mazur's lemma (see \cite{pw} and \cite{yk}), 
  we can deduce    
\begin{theorem}[\cite{jipam}]
Assume that a sequence $\{u_i\}_{i=1}^{+\infty}$ 
 weakly* converges in 
$(L_\mu^\infty({\bf R}^n))^m$ to $u$ as $i$ goes to infinity,  
where $m$ and $n$ are two positive integers. 
Assume that all the values of $u$ and $u_i$ ($i=1,2,3,\cdots$)  
belong to a closed convex set $K$ in ${\bf R}^m$  and that 
$f(w)$ is a continuous convex function from $K$ to  ${\bf R}.$  
Then the inequality (\ref{ineq}) holds 
for any bounded measurable set $\Omega\subset {\bf R}^n.$
\label{th03}
\end{theorem} 

Theorem~\ref{th03} is in fact an extension of  a result given by Ying \cite{yl} 
(or see \cite{jipam} and \cite{jmaa}) and its detail proof can be found in \cite{jipam}.

\vskip0.2cm
\noindent{\small {\bf  Acknowledgement.} 
This work was supported by grants of NSFC 10271121 and 
joint grants of NSFC 10511120278/10611120371 and RFBR 04-02-39026. 
This work was also sponsored by SRF for ROCS, SEM.  
We would like to thank the referee of this paper 
for his/her valuable comments on this work.}
{\small 

}
\end {document}